\documentclass[12pt]{article}
\usepackage{amsmath,amssymb,amsthm}
\usepackage[mathscr]{euscript}
\usepackage[letterpaper]{geometry}
\usepackage{comment}
\usepackage{graphicx}
\usepackage{latexsym}
\usepackage{color}
\usepackage{curves}
\usepackage{psfrag}  

\def\BBox{\kern  -0.16cm\hbox{\vrule width 0.16cm height 0.16cm}}

\renewcommand{\geq}{\geqslant}

\newcommand{\K}{\mathcal{K}} 
 
\newcommand{\E}{\mathbb{E}^3}

\title{Skeletal Geometric Complexes and Their Symmetries}

\author{
Egon Schulte\thanks{Email: schulte@neu.edu}\\[.03in]
Department of Mathematics\\[.03in]
Northeastern University, Boston, MA 02115, USA\\[.08in]
and\\[.08in]
Asia Ivi\'c Weiss\thanks{Email: weiss@mathstat.yorku.ca}\\[.03in]
Department of Mathematics and Statistics\\ 
York University, Toronto, Ontario M3J 1P3, Canada}

\begin{document}
\maketitle

\noindent Polyhedra and polyhedron-like structures have been studied since the early days of geometry (Coxeter~\cite{crp}). The most well-known are the five Platonic solids --- the tetrahedron, octahedron, cube, icosahedron, and dodecahedron. Book XIII of Euclid's ``Elements"~\cite{euclid} was entirely devoted to their mathematical properties. 

\begin{figure}
\centering
$\begin{array}{l}
\includegraphics[width=3in]{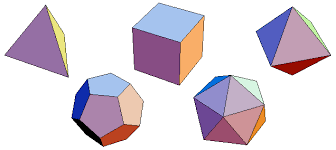}\\[.15in]
\includegraphics[width=.8in]{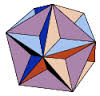}
\includegraphics[width=.8in]{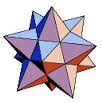}
\includegraphics[width=.8in]{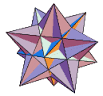}
\includegraphics[width=.8in]{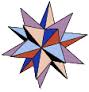}
\end{array}$
\caption{The classical regular polyhedra (Platonic solids and Kepler-Poinsot polyhedra).} 
\label{platokp}
\end{figure}

In modern terminology, the Platonic solids are precisely the regular convex polyhedra in ordinary Euclidean space $\mathbb{E}^3$. Historically, as the name suggests, these figures were viewed as solids bounded by a collection of congruent regular polygons fitting together in a highly symmetric fashion. With the passage of time, the perception of what should constitute a polyhedron has undergone fundamental changes. The famous Kepler-Poinsot star polyhedra no longer are solids but instead form self-intersecting polyhedral structures that feature regular star polygons as faces or vertex-figures. These four ``starry" figures are precisely the regular star polyhedra in $\mathbb{E}^3$~\cite{crp}, which along with the five Platonic solids might be called the classical regular polyhedra~\cite{arp} and are shown in Figure~\ref{platokp}~\cite{wiki,wolfram}.

\begin{figure}
\centering
\includegraphics[width=2.1in]{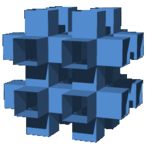}\quad
\includegraphics[width=2.1in]{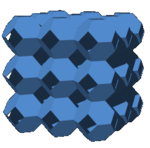}\quad
\includegraphics[width=2.1in]{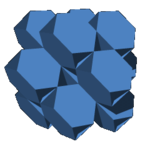}
\caption{The Petrie-Coxeter polyhedra.} 
\label{petcox}
\end{figure}

Why say that a polyhedron must be finite? In the 1930's, Petrie and Coxeter discovered three more regular polyhedra in $\mathbb{E}^3$, each forming a non-compact polyhedral surface which is tessellated by regular convex polygons and bounds two unbounded (congruent) ``solids" (infinite ``handlebodies"). The three polyhedra are shown in Figure~\ref{petcox}~\cite{wiki,wolfram}. The trick here is to allow the vertex-figures to be skew polygons rather than convex polygons. The vertex-figure of a polyhedron $P$ at a vertex $v$ is the polygon whose vertices are the vertices of $P$ adjacent to $v$ and whose edges join two vertices if they lie in a common face of $P$. For example, the faces of the second polyhedron in Figure~\ref{petcox}, denoted by $\{6,4|4\}$, are convex hexagons, four meeting at each vertex; the vertex-figure is a tetragon, a skew square. A little bug moving around a vertex on the polyhedral surface of $\{6,4|4\}$ needs to climb up and down, in an alternating fashion consistent with the fact that the vertex-figure is skew.

\begin{figure}
\centering
\resizebox{180pt}{160pt}{
\begin{picture}(110,100)
\put(0,1){\thicklines\color{red}\line(1,0){75}}
\put(76,0){\thicklines\color{red}\line(0,1){75}}
\put(75,76){\thicklines\color{red}\line(3,2){30}}
\put(105,96){\thicklines\color{red}\line(-1,0){75}}
\put(31,95){\thicklines\color{red}\line(0,-1){75}}
\put(30,21){\thicklines\color{red}\line(-3,-2){30}}
\put(0,-1){\thicklines\color{blue}\line(1,0){75}}
\put(75,1){\thicklines\color{blue}\line(3,2){30}}
\put(106,20){\thicklines\color{blue}\line(0,1){75}}
\put(105,94){\thicklines\color{blue}\line(-1,0){75}}
\put(30,94){\thicklines\color{blue}\line(-3,-2){30}}
\put(1,75){\thicklines\color{blue}\line(0,-1){75}}
\put(0,74){\thicklines\color{green}\line(1,0){75}}
\put(75,74){\thicklines\color{green}\line(3,2){30}}
\put(104,20){\thicklines\color{green}\line(0,1){75}}
\put(105,21){\thicklines\color{green}\line(-1,0){75}}
\put(30,19){\thicklines\color{green}\line(-3,-2){30}}
\put(-1,75){\thicklines\color{green}\line(0,-1){75}}
\put(74,0){\thicklines\color{black}\line(0,1){75}}
\put(0,76){\thicklines\color{black}\line(1,0){75}}
\put(30,96.5){\thicklines\color{black}\line(-3,-2){30}}
\put(29,95){\thicklines\color{black}\line(0,-1){75}}
\put(105,19){\thicklines\color{black}\line(-1,0){75}}
\put(75,-1.5){\thicklines\color{black}\line(3,2){30}}
\multiput(0,0)(75,0){2}{\circle*{4}}
\multiput(0,0)(0,75){2}{\circle*{4}}
\put(75,75){\circle*{4}}
\multiput(30,20)(75,0){2}{\circle*{4}}
\multiput(30,20)(0,75){2}{\circle*{4}}
\put(105,95){\circle*{4}}
\end{picture}}
\caption{The Petrie dual of the cube.} 
\label{petcube}
\end{figure}
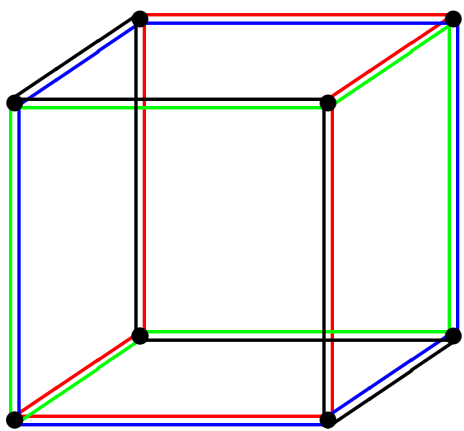

These Petrie-Coxeter polyhedra are precisely the infinite regular polyhedra with convex faces in $\mathbb{E}^3$. We now are at twelve regular polyhedra. Are there more? None were added to the list for another forty years! Then things changed, again! 

In a pioneering paper~\cite{gon}, Gr\"unbaum gave a new definition of polyhedron designed to rid the theory of the psychological block that membranes must be spanning the faces to give a surface. In this skeletal approach, a polyhedron is viewed as a finite or infinite geometric edge graph in $\mathbb{E}^3$ equipped with additional ``polyhedral structure" imposed by a collection of faces, each a polygonal cycle or a two-sided infinite polygonal path. Thus, a polyhedron is comprised of vertices, joined by edges, assembled in careful fashion into polygons, the faces, allowed to be skew or infinite, with the condition that two faces meet at each edge. This approach restored the symmetry in the definition of polyhedron, in that now faces are also allowed to be skew. 

Before moving on to more formal definitions let us consider some examples. 

It turns out that each regular polyhedron $P$ gives rise to one new regular polyhedron~$P^\pi$, the Petrie dual of $P$. The Petrie dual $P^\pi$ has the same edge graph as the original polyhedron $P$ but has as new faces the Petrie polygons of $P$.  A Petrie polygon (or $1$-zigzag) of $P$ is a path along edges such that any two, but no three, consecutive edges belong to a face of $P$. Figure~\ref{petcube} shows the Petrie dual of the ordinary cube, with eight vertices, twelve edges, and four skew hexagonal faces (in red, blue, green, and black), and with three faces meeting at each vertex. Its faces are the four hexagonal Petrie polygons of the cube, each viewed as a polygonal cycle consisting of six edges. Taking Petrie duals of polyhedra is an involutory operation: the Petrie dual of the Petrie dual of $P$ is $P$ itself. 

The nine classical regular polyhedra and their Petrie duals are finite regular polyhedra. In fact, these are the only finite regular polyhedra in $\mathbb{E}^3$. 

Every polyhedron $P$ in $\mathbb{E}^3$ comes with an underlying abstract polyhedron $\mathcal{P}$ describing the incidence among the vertices, edges, and faces of $P$~\cite{arp}. When $P$ has finite faces its abstract polyhedron $\mathcal{P}$ can be realized as a map on a surface. For example, the Petrie dual of the cube shown in Figure~\ref{petcube} is a regular map on the torus, of type $\{6,3\}$, meaning that it has hexagonal faces and triangular vertex-figures. In a map on a surface the faces appear as topological discs, unlike in a geometric polyhedron where by design the membranes are suppressed and the focus is on the edge graph (skeleton).

\begin{figure}
\centering
\includegraphics[width=4.2in]{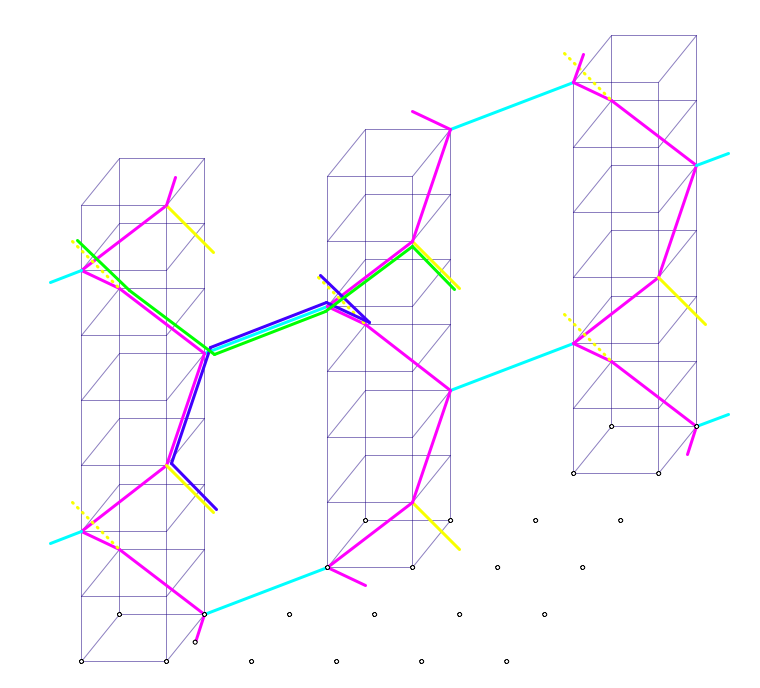}
\caption{The regular polyhedron $\{\infty,3\}^{(b)}$.}
\label{manhattan}
\end{figure}

\begin{figure}
\centering
\includegraphics[width=3.4in]{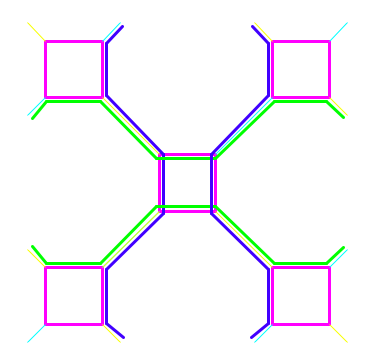}
\caption{Planar projection of $\{\infty,3\}^{(b)}$.}
\label{planarproj}
\end{figure}

Regular polyhedra in $\mathbb{E}^3$ need not have finite faces. For example, the Petrie dual of the Petrie-Coxeter polyhedron $\{6,4|3\}$ is a regular polyhedron whose faces are helical polygons spiraling over a triangular basis, such that four helical faces meet at each vertex. There are several other helix-faced regular polyhedra. Figure~\ref{manhattan} shows an example with square-based helical polygons as faces, three meeting at each vertex. There are three sets of helical polygons:\ ``vertical" helices (in magenta), ``front-to-back" helices (in blue), and ``left-to-right" helices (in green). Figure~\ref{planarproj} illustrates the projection of this polyhedron on a plane perpendicular to the axes of one set of helices, in this case the vertical helices. The vertical helices project onto squares (in magenta), the ``front-to-back" helices to ``top-to-bottom" apeirogons (in blue), and the ``left-to-right" helices onto ``left-to-right" apeirogons (in green).

Loosely speaking there are forty-eight regular polyhedra in $\mathbb{E}^3$, up to similarity, and these are known as the Gr\"unbaum-Dress polyhedra. More precisely, there are thirty-six individual polyhedra as well as twelve one-parameter families in which the polyhedra look very much alike. Gr\"unbaum~\cite{gon} found all polyhedra, save one, and Dress~\cite{dr1,dr2} discovered the final instance and proved the completeness of the enumeration; a quick method of arriving at the full characterization is described by McMullen and Schulte in~\cite[Sect. 7E]{arp} (and \cite{rpo}).

While regular polyhedra are all (maximally) symmetric by reflection, chiral polyhedra only have maximum symmetry by rotation (in a sense described below). It appears that the term ``chiral" was first used by Lord Kelvin in 1893~\cite{kelvin}, who wrote, ``I call any geometrical figure, or group of points, chiral, and say that it has chirality if its image in a plane mirror, ideally realized, cannot be brought to coincide with itself". The word ``chiral" comes from the Greek word $\chi\varepsilon\iota\rho$ (kheir), which means ``hand", in recognition of the fact that human hands are chiral. Many geometric figures and many objects in the sciences exhibit chirality. Our use of the term ``chiral" is more restrictive and applies to structures that have chirality features in the presence of otherwise maximum possible symmetry. A snub cube, a perfect example of a polyhedron which is chiral according to Lord Kelvin's definition, will not be a chiral polyhedron according to our definition as it fails to have maximum possible symmetry. For example, it has faces of two types, namely squares and equilateral triangles.

That chiral polyhedra really exist is by no means obvious. In fact, there are no classical examples of chiral polyhedra which could serve as inspiration to discover new examples. Nevertheless, there are very many chiral polyhedra in $\mathbb{E}^3$. They can be completely classified into six large one-parameter families, up to similarity \cite{chir1,chir2}. 
These polyhedra have stunning properties, some quite counter-intuitive. 

\begin{figure}
\centering
\includegraphics[width=3.6in]{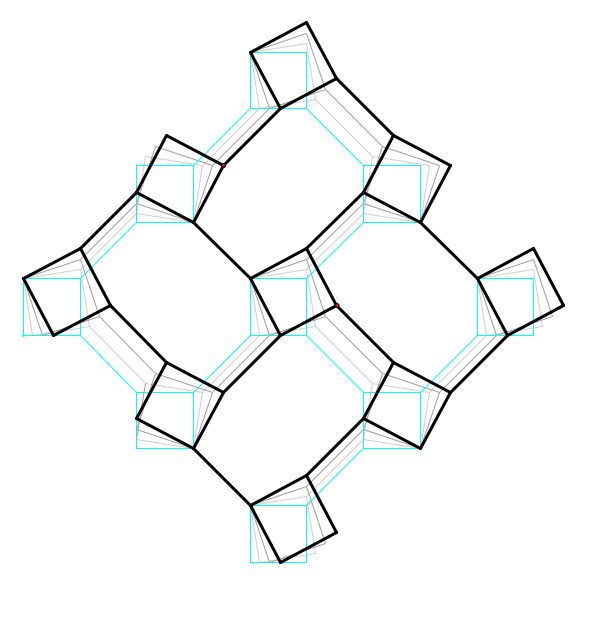}
\caption{Deforming a chiral helix-faced polyhedron.}
\label{chirdeform}
\end{figure}

It turns out that there are three large families, each with one rational parameter (when described up to similarity), in which the faces and vertex-figures are finite skew polygons; in each family, essentially any two polyhedra are combinatorially non-isomorphic. There are also three large  families, each with one real parameter (when described up to similarity), in which the faces are helical polygons and the vertex-figures are convex polygons; now, in each family, any two polyhedra are combinatorially isomorphic, and each polyhedron is isomorphic to a regular helix-faced polyhedron into which it can be deformed continuously \cite{pelwei}. Figure~\ref{chirdeform} below employs a planar projection of a chiral helix-faced polyhedron to illustrate how it can be moved, in fact continuously, into the regular helix-faced polyhedron shown in Figure~\ref{manhattan}; the chiral polyhedra used in this process all have only half as many symmetries as the regular polyhedron.

In a polyhedron, only two faces are permitted to meet at an edge. There is a broader type of polyhedral structure in $\mathbb{E}^3$, called polygonal complex, in which more than two faces are allowed to meet at an edge. Like a polyhedron, a polygonal complex $K$ is comprised of vertices, edges, and polygonal faces. However, unlike in a polyhedron, the vertex-figures may be finite graphs which are not cycles; the vertex-figure at a vertex $v$ of $K$ is the graph whose vertices are the neighbors of $v$ in $K$ and whose edges join two vertices if they lie in a common face of $K$. The polyhedra are precisely the polygonal complexes in which exactly two faces meet at each edge. A familiar example is the complex consisting of all square faces of the standard cubical tessellation of $\mathbb{E}^3$ (see Figure~\ref{2skeleton} below); its faces are squares and its vertex-figures are edge graphs of octahedra. 

\begin{figure}
\centering
\resizebox{260pt}{192pt}{
\begin{picture}(140,110)
\multiput(10,0)(0,40){3}{
\begin{picture}(100,25)
\thinlines
\multiput(0,0)(40,0){3}{\circle*{3}}
\multiput(18,12)(40,0){3}{\circle*{3}}
\multiput(36,24)(40,0){3}{\circle*{3}}
\multiput(0,0)(18,12){3}{\line(1,0){80}}
\multiput(0,0)(40,0){3}{\line(3,2){36}}
\end{picture}}
\put(10,0){
\begin{picture}(180,140)
\thinlines
\multiput(0,0)(40,0){3}{\line(0,1){80}}
\multiput(18,12)(40,0){3}{\line(0,1){80}}
\multiput(36,24)(40,0){3}{\line(0,1){80}}
\end{picture}}
\put(10,0){
\begin{picture}(100,80)
\thicklines
\put(40,-1){\color{green}\line(1,0){40}}
\multiput(80,0)(.45,.3){40}{\color{green}\circle*{1}}          
\put(98,12){\color{green}\line(0,1){40}}
\put(98,51){\color{green}\line(-1,0){40}}
\multiput(58,50.5)(-.45,-.3){40}{\color{green}\circle*{1}}
\put(40,40){\color{green}\line(0,-1){40}}
\end{picture}}
\put(10,0){
\begin{picture}(100,80)
\thicklines
\put(58,49){\color{blue}\line(1,0){40}} 
\multiput(98,52)(.45,.3){40}{\color{blue}\circle*{1}}          
\put(116,64){\color{blue}\line(0,1){40}}
\put(116,104){\color{blue}\line(-1,0){40}}
\multiput(76,104)(-.45,-.3){40}{\color{blue}\circle*{1}}
\put(61,92){\color{blue}\line(0,-1){40}}
\end{picture}}
\put(10,0){
\begin{picture}(100,80)
\thicklines
\put(40,1){\color{yellow}\line(1,0){40}}
\put(80,0){\color{yellow}\line(0,1){40}}
\multiput(80,40)(.45,.3){40}{\color{yellow}\circle*{1}}    
\put(98,53.5){\color{yellow}\line(-1,0){40}}
\put(56.5,52){\color{yellow}\line(0,-1){40}}
\multiput(40,0)(.45,.3){40}{\color{yellow}\circle*{1}}
\end{picture}}
\put(10,0){
\begin{picture}(100,80)
\thicklines
\put(58,55){\color{red}\line(1,0){40}}
\put(98,52){\color{red}\line(0,1){40}}
\multiput(98,92)(.45,.3){40}{\color{red}\circle*{1}}    
\put(116,105){\color{red}\line(-1,0){40}}
\put(76,104){\color{red}\line(0,-1){40}}
\multiput(58,55)(.45,.3){40}{\color{red}\circle*{1}}
\end{picture}}
\put(10,0){
\begin{picture}(100,80)
\thicklines
\put(56.5,52){\color{yellow}\line(0,1){40}}
\multiput(58,53.5)(.45,.3){40}{\color{yellow}\circle*{1}}  
\put(58,92){\color{yellow}\line(1,0){40}}
\multiput(98,94)(.45,.3){40}{\color{yellow}\circle*{1}}
\put(76,64){\color{yellow}\line(1,0){40}}
\put(117,64){\color{yellow}\line(0,1){40}}
\put(61,52){\color{blue}\line(0,-1){40}}
\multiput(58,49)(-.45,-.3){40}{\color{blue}\circle*{1}}
\put(40,40){\color{blue}\line(1,0){40}}
\put(79,40){\color{blue}\line(0,-1){40}}
\multiput(80,-1)(.45,.3){40}{\color{blue}\circle*{1}}  
\put(98,12){\color{blue}\line(-1,0){40}}
\put(58,55){\color{red}\line(-1,0){40}} 
\multiput(18,51)(.45,.3){40}{\color{red}\circle*{1}}
\multiput(58,13.5)(.45,.3){40}{\color{red}\circle*{1}}
\put(36,23){\color{red}\line(1,0){40}}
\put(37,24){\color{red}\line(0,1){40}}
\put(55,52){\color{red}\line(0,1){40}}
\put(58,49){\color{blue}\line(-1,0){40}} 
\put(18,12){\color{blue}\line(0,1){40}}
\multiput(58,49)(.45,.3){40}{\color{blue}\circle*{1}}
\put(76,64){\color{blue}\line(0,-1){40}}
\put(76,25){\color{blue}\line(-1,0){40}}
\multiput(36,24)(-.45,-.3){40}{\color{blue}\circle*{1}}
\multiput(58,50.5)(.45,.3){40}{\color{green}\circle*{1}}
\multiput(18,13)(.45,.3){40}{\color{green}\circle*{1}}
\put(18,12){\color{green}\line(1,0){40}}
\put(59.5,12){\color{green}\line(0,1){40}}
\put(76,64){\color{green}\line(-1,0){40}}
\put(36,64){\color{green}\line(0,-1){40}}
\multiput(58,53.5)(-.45,-.3){40}{\color{yellow}\circle*{1}}
\put(58,53.5){\color{yellow}\line(-1,0){40}}
\put(40,40){\color{yellow}\line(0,1){40}}
\put(40,80){\color{yellow}\line(-1,0){40}}
\multiput(0,80)(.45,.3){40}{\color{yellow}\circle*{1}}
\put(18,92){\color{yellow}\line(0,-1){40}}
\put(40,40){\color{red}\line(-1,0){40}}
\put(1,40){\color{red}\line(0,1){40}}
\multiput(40,42)(.45,.3){40}{\color{red}\circle*{1}}
\put(58,92){\color{red}\line(-1,0){40}}
\multiput(18,93)(-.45,-.3){40}{\color{red}\circle*{1}}
\put(55,53){\color{red}\line(0,-1){40}}
\put(0,39){\color{green}\line(0,1){40}}
\put(0,81){\color{green}\line(1,0){40}}
\multiput(40,80)(.45,.3){40}{\color{green}\circle*{1}}
\multiput(0,40)(.45,.3){40}{\color{green}\circle*{1}}
\put(59.5,92){\color{green}\line(0,-1){40}}
\put(58,50.5){\color{green}\line(-1,0){40}}
\put(58,52){\circle*{6}}
\put(18,52){\circle*{6}}
\put(98,52){\circle*{6}}
\put(58,12){\circle*{6}}
\put(58,92){\circle*{6}}
\put(40,40){\circle*{6}}
\put(76,64){\circle*{6}}
\end{picture}}
\multiput(10,0)(0,40){3}{
\begin{picture}(100,25)
\thinlines
\multiput(0,0)(40,0){3}{\circle*{3}}
\multiput(18,12)(40,0){3}{\circle*{3}}
\multiput(36,24)(40,0){3}{\circle*{3}}
\multiput(0,0)(18,12){3}{\line(1,0){80}}
\multiput(0,0)(40,0){3}{\line(3,2){36}}
\end{picture}}
\end{picture}}
\caption{A local picture of the polygonal complex $K_{4}(1,2)$, with skew hexagonal faces.}
\label{figk412}
\end{figure}
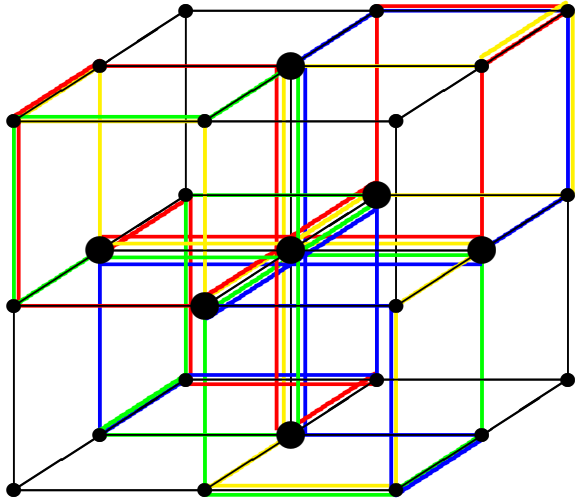

A more complicated example is illustrated in Figure~\ref{figk412}; the faces are the (hexagonal) Petrie polygons of alternate cubes in the standard cubical tessellation, and four faces meet at each edge. Shown are the twelve faces that have a vertex in common, here located at the center. This is the regular polygonal $K_{4}(1,2)$ described later in the text.

In addition to the forty-eight regular polyhedra in $\mathbb{E}^3$ there are twenty-five regular polygonal complexes, all infinite; they were described by Pellicer and Schulte in~\cite{pesch1,pesch2}. Their edge graphs are nets well-known to crystallographers~\cite{okeef,schnets}.

The present article focuses on skeletal polyhedral structures in ordinary space $\E$. For important results about skeletal structures in higher dimensional Euclidean spaces we refer the reader to \cite{ar,bra,bhp,pm} and in particular to the forthcoming monograph by McMullen~\cite{mr}. 
\bigskip

\noindent
{\bf Polygons\/}
\medskip

\noindent
A study of polyhedra must begin with polygons, as these occur as faces and vertex-figures. Although in the present context only regular polygons occur, it is useful to give a formal definition of polygon (see \cite{gon}). 

A finite polygon, or $n$-gon (with $n\geq 3$), consists of a sequence $(v_1, v_2, \dots, v_n)$ of $n$ distinct points in $\mathbb{E}^3$, together with the line segments $(v_1, v_2), (v_2,v_3), \ldots, (v_{n-1},v_n)$ and $(v_n, v_1)$. Similarly, a (discrete) infinite polygon, or apeirogon, consists of an infinite sequence of distinct points $(\dots, v_{-2},v_{-1}, v_0, v_1, v_2, \dots)$ in $\mathbb{E}^3$, along with the line segments $(v_i, v_{i+1})$ for each integer $i$, such that each compact subset of $\mathbb{E}^3$ meets only finitely many line segments. In either case the points are the vertices and the line segments the edges of the polygon. It is important that the points be distinct. Topologically, a finite polygon is a 1-sphere and an infinite polygon is a real line.

A polygon is regular if its symmetry group acts transitively on the set of flags (incident vertex-edge pairs). Recall that the symmetry group of a geometric figure in $\mathbb{E}^3$ consists of the isometries of $\mathbb{E}^3$ which leave the figure as a whole invariant. Regular polygons in $\mathbb{E}^3$ are either finite, planar (convex or star-) polygons or non-planar (skew) polygons, or infinite, planar zigzags or helical polygons. For example, the faces of the polyhedron in  
Figure~\ref{petcube} are skew hexagons. The vertices of a helical polygon lie on a helix, and those of a zigzag polygon lie alternately on two parallel lines (if the lines coincide, the polygon becomes linear).
\bigskip

\noindent
{\bf Polyhedra\/}
\medskip

\noindent
Following Gr\"unbaum as in~\cite{gon}, a polyhedron $P$ in $\mathbb{E}^3$ consists of a set $V$ of points, called vertices, a set $E$ of line segments, called edges, and a set $F$ of polygons, called faces, satisfying the following properties:\\[.04in]
(a)\ The graph $(V,E)$, called the edge graph of $P$, is connected.\\[.01in]
(b)\ The vertex-figure at each vertex of $P$ is connected. By the vertex-figure of $P$ at a vertex~$v$ we mean the graph whose vertices are the neighbors of $v$ in the edge graph of~$P$ and whose edges are the line segments $(u,w)$, where $(u,v)$ and $(v,w)$ are edges of a common face of $P$.\\[.01in]
(c)\ Each edge of $P$ is contained in exactly two faces of $P$.\\[.01in]
(d)\ $P$ is discrete, meaning that each compact subset of $\mathbb{E}^3$ meets only finitely many faces of~$P$.
\smallskip

A flag of a polyhedron $P$ is an incident vertex-edge-face triple. Two flags are called adjacent if they differ in precisely one element, that is, in the vertex, the edge, or the face; accordingly, the flags are also called $0$-adjacent, $1$-adjacent, or $2$-adjacent. 

A polyhedron $P$ is regular if its symmetry group $G(P)$ acts transitively on the set of flags. A polyhedron $P$ is chiral if $G(P)$ has exactly two orbits on the flags such that adjacent flags are always in distinct orbits.

If $P$ is regular or chiral, then $G(P)$ acts transitively, separately, on the sets of faces, edges, and vertices of $P$, so in particular $P$ has mutually congruent faces, edges, and vertex-figures, respectively. The faces of $P$ are finite or infinite regular polygons, and the vertex-figures are finite regular polygons (by discreteness). Each regular or chiral polyhedron $P$ is assigned a basic Schl\"afli symbol, or type, $\{p,q\}$, where $p$ is the number of vertices in a face (and possibly is $\infty$) and $q$ is the (finite) number of faces meeting at a vertex. 
\bigskip

\noindent
{\bf Regular Polyhedra\/}
\medskip

\noindent
The symmetry group $G(P)$ of a regular polyhedron $P$ is generated by three reflections $R_0,R_1,R_2$ in points, lines, or planes. (Reflections in lines are half-turns.) These reflections are determined by the way in which they act on a chosen base flag of $P$. More specifically, $R_0$, $R_1$, and $R_2$ map this base flag to its $0$-adjacent, $1$-adjacent, or $2$-adjacent flag, respectively. For example, for the Petrie dual of the cube shown in Figure~\ref{petcube}, $R_0$ is a half-turn and $R_1$ and $R_2$ are plane reflections. Only in the case of the nine classical regular polyhedra are all three generators plane reflections. If $P$ is of type $\{p,q\}$, then $G(P)$ satisfies the Coxeter relations 
\begin{equation}
\label{coxrelations}
R_{0}^{2}=R_{1}^{2}=R_{2}^{2}=(R_{0}R_{1})^{p}
=(R_{1}R_{2})^{q} = (R_{0}R_{2})^{2} = 1 ,
\end{equation}
but, in most case, other independent relations, too. 

It is customary to designate these polyhedra by generalized Schl\"afli symbols that usually are obtained by padding the basic symbol $\{p,q\}$ with additional symbols signifying specific information (such as extra defining relations for the symmetry group). For example, the symbol for the Petrie dual of the cube is $\{6,3\}_4$, with the subscript $4$ attached to the basic symbol $\{6,3\}$ indicating that the symmetry group has a presentation consisting of the Coxeter relations in (\ref{coxrelations}) and the single extra relation $(R_{0}R_{1}R_{2})^{4}=1$; the symmetry $R_{0}R_{1}R_{2}$ shifts a certain Petrie polygon of $\{6,3\}_4$ one step along itself, forcing the Petrie polygon to close up after four steps and hence have length $4$.

Following the classification scheme described in \cite{rpo,arp}, the forty-eight regular polyhedra in $\E$ broadly fall into certain families  as follows. There are eighteen finite polyhedra and thirty apeirohedra (infinite polyhedra). The eighteen finite polyhedra are the nine classical regular polyhedra and their Petrie duals. The thirty apeirohedra further fall into three families:\ the six planar polyhedra, the twelve ``reducible" (blended) polyhedra, and the twelve ``irreducible" (pure) polyhedra. 

\begin{figure}
\centering
\includegraphics[width=3in]{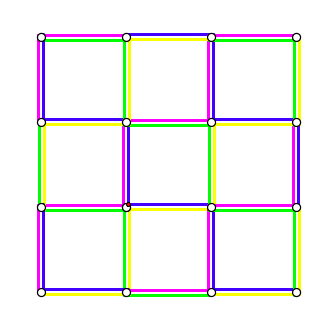}
\caption{The polyhedron $\{\infty,4\}_4$}
\label{inftyfour}
\end{figure}

The six planar polyhedra are the familiar regular plane tessellations $\{3,6\}$, $\{6,3\}$, and $\{4,4\}$, by triangles, hexagons, and squares, respectively, and their Petrie duals $\{\infty,6\}_3$, $\{\infty,3\}_6$, and $\{\infty,4\}_4$. The polyhedron $\{\infty,4\}_4$ is shown in Figure~\ref{inftyfour}; it has the same vertices and edges as the square tessellation $\{4,4\}$, but its faces are the Petrie polygons of $\{4,4\}$, which are zigzag polygons, four meeting at each vertex.

The ``reducible" polyhedra are precisely the regular polyhedra $P$ whose symmetry groups $G(P)$ acts (affinely) reducibly on $\mathbb{E}^3$, meaning that there is a non-trivial linear subspace $L$ (a line or plane) which is invariant in the sense that $G(P)$ permutes the translates of $L$. 

The twelve reducible polyhedra are blends, in the sense that they are obtained by ``blending" one of the six planar polyhedra with one of the two linear regular polygons, that is, with either a line segment $\{\;\}$ or a linear apeirogon $\{\infty\}$, contained in a line perpendicular to the plane with the planar polyhedron. The projections of a blended polyhedron onto its two component subspaces recover the two original components, that is, the original planar polyhedron as well as the line segment or apeirogon. A blend is determined up to similarity by the ratio of the edge lengths in the two components. Thus each blended regular polyhedron really represents a one-parameter family of mutually non-similar regular polyhedra, which all look very much alike, in the sense that they have the same combinatorial characteristics. 

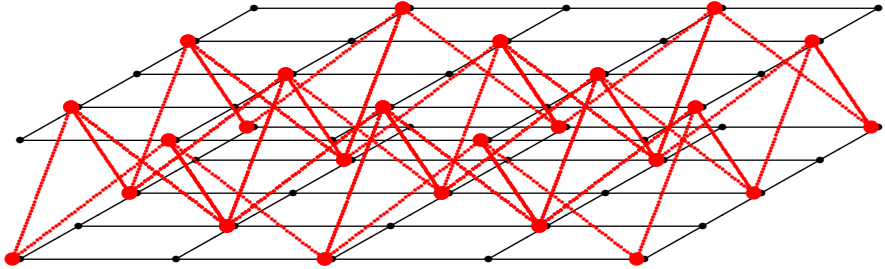
\begin{figure}
\centering
\resizebox{360pt}{100pt}{
\begin{picture}(240,80)
\multiput(0,0)(0,36){2}{
\begin{picture}(240,50)
\thinlines
\multiput(0,0)(40,0){5}{\circle*{2}}
\multiput(15,10)(40,0){5}{\circle*{2}}
\multiput(30,20)(40,0){5}{\circle*{2}}
\multiput(45,30)(40,0){5}{\circle*{2}}
\multiput(60,40)(40,0){5}{\circle*{2}}
\multiput(0,0)(15,10){5}{\line(1,0){160}}
\multiput(0,0)(40,0){5}{\line(3,2){60}}
\end{picture}}
\multiput(192,20)(0.3,0.92){50}{\color{red}\circle*{1}}
\multiput(0,0)(30,20){2}{
\multiput(0,0)(80,00){2}{
\multiput(2,0)(0.3,0.92){50}{\color{red}\circle*{1}}
\multiput(2,0)(0.8,0.72){50}{\color{red}\circle*{1}}
\multiput(17,46)(0.8,-0.72){50}{\color{red}\circle*{1}} 
\multiput(42,36)(0.3,-0.52){50}{\color{red}\circle*{1}} 
\put(2,0){\color{red}\circle*{4}}
\put(17,46){\color{red}\circle*{4}}
\put(42,36){\color{red}\circle*{4}}
\put(57,10){\color{red}\circle*{4}}}}
\multiput(40,0)(80,0){2}{
\multiput(2,36)(0.8,-0.72){50}{\color{red}\circle*{1}} 
\multiput(2,36)(0.3,-0.52){50}{\color{red}\circle*{1}} 
\multiput(42,0)(0.3,0.92){50}{\color{red}\circle*{1}} 
\multiput(17,10)(0.8,0.72){50}{\color{red}\circle*{1}} 
\put(2,36){\color{red}\circle*{4}}
\put(57,46){\color{red}\circle*{4}}
\put(42,0){\color{red}\circle*{4}}
\put(17,10){\color{red}\circle*{4}}}         
\multiput(15,10)(30,20){2}{
\multiput(40,0)(80,00){2}{
\multiput(2,0)(0.3,0.92){50}{\color{red}\circle*{1}}
\multiput(2,0)(0.8,0.72){50}{\color{red}\circle*{1}}
\multiput(17,46)(0.8,-0.72){50}{\color{red}\circle*{1}} 
\multiput(42,36)(0.3,-0.52){50}{\color{red}\circle*{1}} 
\put(2,0){\color{red}\circle*{4}}
\put(17,46){\color{red}\circle*{4}}
\put(42,36){\color{red}\circle*{4}}
\put(57,10){\color{red}\circle*{4}}}
\multiput(0,0)(80,0){2}{
\multiput(2,36)(0.8,-0.72){50}{\color{red}\circle*{1}} 
\multiput(2,36)(0.3,-0.52){50}{\color{red}\circle*{1}} 
\multiput(42,0)(0.3,0.92){50}{\color{red}\circle*{1}} 
\multiput(17,10)(0.8,0.72){50}{\color{red}\circle*{1}} 
\put(2,36){\color{red}\circle*{4}}
\put(57,46){\color{red}\circle*{4}}
\put(42,0){\color{red}\circle*{4}}
\put(17,10){\color{red}\circle*{4}}}}
\end{picture}}
\caption{The blend of the square tessellation with the line segment.}
\label{blendline}
\end{figure}

Figure~\ref{blendline} illustrates the blend of the square tessellation $\{4,4\}$ in the $xy$-plane of $\mathbb{E}^3$ with a line segment $[-1,1]$ (say) positioned along the $z$-axis in $\mathbb{E}^3$. The vertices of the polyhedron lie in the planes $z=-1$ and $z=1$, and are obtained from those of $\{4,4\}$ by alternately raising or lowering vertices. Its faces are tetragons (skew squares, one over each square of the original tessellation), with vertices alternating between the two planes. Its designation is $\{4,4\}\#\{\;\}$, where $\#$ indicates the blending operation. 

The blend of the square tessellation in the $xy$-plane with a linear apeirogon positioned along the $z$-axis has helical polygons as faces; these rise in two-sided infinite vertical towers above the squares of the original  tessellation, such that the helical polygons over adjacent squares have opposite orientations (left-handed or right-handed) and meet along every fourth edge. The designation in this case is $\{4,4\}\#\{\infty\}$. 

The ``irreducible", or pure, regular polyhedra (now apeirohedra) have symmetry groups $G(P)$ that act (affinely) irreducibly on $\mathbb{E}^3$. There are exactly twelve pure polyhedra, all afforded by the rich geometry of the standard cubical tessellation in $\mathbb{E}^3$. They are related to each other by a number of geometric operations such as duality, Petrie-duality, and other such operations not further discussed here.  

The most prominent pure polyhedra are the three Petrie-Coxeter polyhedra $\{4,6\,|\,4\}$, $\{6,4\,|\,4\}$ and $\{6,6\,|\,3\}$ shown in Figure~\ref{petcox}. The first two entries in a symbol record the standard Schl\"afli symbol, while the last entry gives the length $h$ of the holes (formed by edge paths which leave a vertex by the second edge on the left from which they entered), which in this case is $4$ or $3$. The symmetry group now has a presentation consisting of the Coxeter relations in (\ref{coxrelations}) and the single extra relation $(R_{0}R_{1}R_{2}R_{1})^{h}=1$, with $h=4$ or $3$; the symmetry $R_{0}R_{1}R_{2}R_{1}$ shifts the edges of a certain hole one step along the hole, forcing it to close up after $h$ steps. These three polyhedra are the only infinite regular polyhedra in $\E$ with convex faces. 

\begin{figure}
\centering
\includegraphics[width=3.5in]{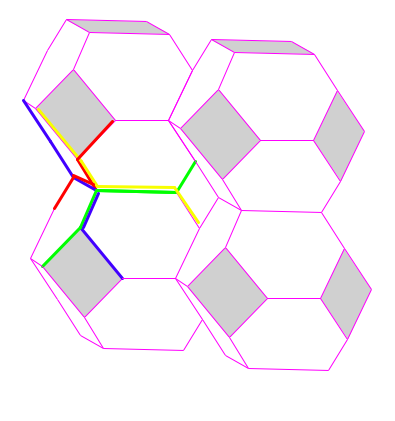}\\[-.45in]
\caption{The Petrie dual of $\{6,4\,|\,4\}$.}
\label{pet643}
\end{figure}

The Petrie duals have helical faces given by the Petrie polygons of the original polyhedron. The Petrie dual of $\{6,4\,|\,4\}$ is depicted in Figure~\ref{pet643} and is designated by $\{\infty,4\}_{6,4}$. Its symmetry group has a presentation given by the relations in (\ref{coxrelations}) and the two extra relations $(R_{0}R_{1}R_{2})^{6}=1$ and $(R_{0}(R_{1}R_{2})^{2})^{4}=1$; the exponents $6$ and $4$, respectively, are the lengths of the Petrie polygons or $2$-zigzags (edge-paths which leave a vertex from the second edge from which they entered, but in the oppositely oriented sense at alternate vertices). 

It is an intriguing question why there are only twelve pure regular apeirohedra. The answer is that there are four possible mirror vectors for the generators, and three ways of grouping the polyhedra by crystallographic Platonic solids that are associated with them, as further explained below. Table~\ref{break} presents the breakdown by mirror vector and Platonic solid, and as well gives information about the nature of the faces and vertex-figures. 

If $P$ is a regular polyhedron and $R_0,R_1,R_2$ are the generators for $G(P)$, the mirror vector $(d_0,d_1,d_2)$ of $P$ records the dimensions $d_0$, $d_1$ and $d_2$ of the mirrors (fixed point sets) of $R_0$, $R_1$ and $R_2$, respectively. Only the four mirror vectors $(2,1,2)$, $(1,1,2)$, $(1,2,1)$ and $(1,1,1)$ can occur, mostly due to the irreducibility of $G(P)$. For example, for the three polyhedra with mirror vector $(1,1,1)$ the symmetry group is generated by three half-turns and therefore contains only proper isometries; thus these helix-faced polyhedra occur geometrically in two enantiomorphic forms, yet they are geometrically regular, not chiral! The polyhedron $\{\infty,3\}^{(b)}$ was shown in Figure~\ref{manhattan}.

To explain the grouping by Platonic solid, suppose again that $R_0,R_1,R_2$ are the distinguished generators of $G(P)$. Then $R_1$ and $R_2$, but not $R_0$, fix the vertex in the base flag, the origin $o$ (say), of $P$. Translate the mirror of $R_0$ so that it passes through $o$, and denote the reflection in the translated mirror by $R_{0}'$. Then $G':=\langle R_0',R_1,R_2\rangle$ is a finite irreducible group of isometries isomorphic to the special group of $G(P)$, the quotient of $G$ by its translation subgroup. Now alter (if needed) the generators $R_0',R_1,R_2$ as follows. If a generator is a half-turn (with $1$-dimensional mirror), replace it by the reflection in the plane through $o$ perpendicular to its rotation axis; otherwise leave the generator unchanged. The resulting plane reflections $\widehat{R}_{0},\widehat{R}_{1},\widehat{R}_{2}$, respectively, generate a finite irreducible reflection group $\widehat{G}$ in $\E$, which is crystallographic by the discreteness of~$P$. There are only three possible choices for $\widehat{G}$ and its generators, namely the symmetry groups of the tetrahedron $\{3,3\}$, the octahedron $\{3,4\}$, or the cube $\{4,3\}$. In this way, each pure regular apeirohedron is associated with $\{3,3\}$, $\{3,4\}$, or $\{4,3\}$. 

\begin{table}
\centering
\begin{tabular}{|c||c|c|c|c|c|} 
\hline     
 Mirror    &  $\{3,3\}$ & $\{3,4\}$ & $\{4,3\}$  &faces&vertex-\\
 vector    &                &                &                &         &figures\\
\hline\hline
(2,1,2) & $\{6,6 | 3\}$ & $\{6,4 | 4\}$ & $\{4,6 | 4\}$ &planar & skew  \\
(1,1,2) & $\{\infty,6\}_{4,4}$ & $\{\infty,4\}_{6,4}$ & $\{\infty,6\}_{6,3}$ &helical & skew\\
(1,2,1) & $\{6,6\}_{4}$ & $\{6,4\}_{6}$ & $\{4,6\}_{6}$&  skew &planar \\
(1,1,1) & $\{\infty,3\}^{(a)}$ & $\{\infty,4\}_{\cdot,*3}$ &
$\{\infty,3\}^{(b)}$ &helical &planar  \\
\hline
\end{tabular}
\caption{The pure regular polyhedra.}
\label{break}
\end{table}
\bigskip

\eject
\noindent
{\bf Chiral Polyhedra\/}
\medskip

\noindent
Chirality is a fascinating phenomenon that does not occur in the classical theory of polyhedra and polytopes. According to the definition, a polyhedron $P$ in $\E$ is chiral if its symmetry group $G(P)$ has two orbits on the flags such that any two adjacent flags are in distinct orbits. This ``global" definition also has a ``local" counterpart:\ a polyhedron $P$ is chiral if and only if $P$ is not regular but, for some incident vertex-face pair of $P$ (the base vertex-face pair), there exist two symmetries $S_{1}$ and $S_{2}$ of $P$ such that $S_1$ fixes the polygonal face in this pair and cyclically permutes the vertices of that face, whereas $S_2$ fixes the vertex in this pair and cyclically permutes the vertices of the vertex-figure of $P$ at that vertex. (Then such symmetries exist for every incident vertex-face pair.) By replacing one of these symmetries by its inverse, if need be, we can further assume that $T:=S_{1}S_{2}$ is a symmetry that interchanges the two vertices of an edge (the base edge) and the two faces meeting at that edge; combinatorially, $T$ acts like a half-turn about the center of that edge. 

Thus a non-regular polyhedron $P$ is chiral if there are two geometric symmetries $S_1$ and $S_2$ of $P$ of the above kind that act combinatorially like rotations about a face and an incident vertex, respectively. However, $S_1$ and $S_2$ are usually not geometric rotations. In fact, each of these symmetries must be an ordinary rotation, a rotatory reflection (a rotation followed by a reflection in a plane perpendicular to the rotation axis), or a screw motion (a rotation followed by a translation along the rotation axis). 

Note that a regular polyhedron also has two symmetries like $S_1$ and $S_2$. In fact, if $R_0,R_1,R_2$ are the distinguished generators of $G(P)$, we may take $S_{1}:=R_0R_1$ and $S_{2}:=R_1R_2$ and observe that then $T=R_0R_2$.  In this case the symmetries $S_1$ and $S_2$ generate the (combinatorial) rotation subgroup $G^+(P)$ of $G(P)$, which usually does not only consist of geometric rotations.

The symmetry group $G(P)$ of a chiral polyhedron $P$ is generated by the two symmetries $S_1$ and $S_2$. If $P$ is of type $\{p,q\}$, then $G(P)$ satisfies the standard relations (for the rotation groups of Coxeter groups),
\begin{equation}
\label{rotcoxrels}
S_{1}^{p}=S_{2}^{q}=(S_{1}S_{2})^{2} = 1 ,
\end{equation}
and in general other independent relations too. 

The classification of chiral polyhedra is rather involved~\cite{chir1,chir2}. There are no finite examples, no planar examples, and also no blended examples. Thus the symmetry group of a chiral polyhedron is an affinely irreducible crystallographic group in $\E$. 

The chiral polyhedra in $\E$ fall into six infinite families, each with one or two free parameters depending on whether the classification is up to similarity or congruence, respectively. The six families comprise three families of polyhedra with finite skew faces and three families of polyhedra with infinite helical faces. It is helpful to slightly extend each family by allowing the parameters to take certain normally excluded values which would make the corresponding polyhedron regular, not chiral, and in some cases also finite. Each extended family will then contain exactly two regular polyhedra (possibly, one finite polyhedron), while all other polyhedra are chiral polyhedra.

\begin{figure}
\centering
\resizebox{230pt}{233pt}{
\begin{picture}(260,265)
\multiput(0,10)(0,90){3}{
\begin{picture}(180,45)
\thinlines
\multiput(0,0)(90,0){3}{\circle*{6}}
\multiput(37.5,25)(90,0){3}{\circle*{6}}
\multiput(75,50)(90,0){3}{\circle*{6}}
\multiput(0,0)(37.5,25){3}{\line(1,0){180}}
\multiput(0,0)(90,0){3}{\line(3,2){75}}
\end{picture}}
\put(0,10){
\begin{picture}(270,225)
\put(127.5,115){\circle*{10}}
\end{picture}}
\put(0,10){
\begin{picture}(270,225)
\thinlines
\multiput(0,0)(90,0){3}{\line(0,1){180}}
\multiput(37.5,25)(90,0){3}{\line(0,1){180}}
\multiput(75,50)(90,0){3}{\line(0,1){180}}
\end{picture}}
\put(0,10){   
\begin{picture}(270,225)
\thicklines
\multiput(-0.5,0)(1,0){3}{\color{red}\line(0,1){90}}
\multiput(-0.5,0)(1.4,0){3}{\color{red}\line(3,2){37.5}}
\multiput(89,88)(1.4,0){3}{\color{red}\line(3,2){37.5}}
\multiput(125.5,25)(1,0){3}{\color{red}\line(0,1){90}}
\multiput(0,89.5)(0,1){3}{\color{red}\line(1,0){90}}
\multiput(37.5,24.5)(0,1){3}{\color{red}\line(1,0){90}}
\end{picture}}
\put(0,10){ 
\begin{picture}(270,225)
\multiput(90,90)(0,1.5){60}{\color{green}\circle*{3}}
\multiput(0,180)(1.5,0){60}{\color{green}\circle*{3}}
\multiput(0,180)(0.625,0.415){60}{\color{green}\circle*{3}}
\multiput(37.5,115)(0,1.5){60}{\color{green}\circle*{3}}
\multiput(90,92)(0.625,0.415){60}{\color{green}\circle*{3}}
\multiput(37.5,117)(1.5,0){60}{\color{green}\circle*{3}}
\end{picture}}
\put(0,10){ 
\begin{picture}(270,225)
\multiput(75,140)(0,1.5){60}{\color{blue}{\circle*{3}}}
\multiput(75,230)(1.5,0){60}{\color{blue}{\circle*{3}}}
\multiput(37.5,115)(0.625,0.415){60}{\color{blue}{\circle*{3}}}
\multiput(127.5,205)(0.625,0.415){60}{\color{blue}{\circle*{3}}}
\multiput(37.5,113)(1.5,0){60}{\color{blue}{\circle*{3}}}
\multiput(125.5,115)(0,1.5){60}{\color{blue}{\circle*{3}}}
\end{picture}}
\put(255,240){   
\begin{picture}(270,225)
\thicklines\thicklines
\multiput(0.5,0)(-1,0){3}{\color{red}\line(0,-1){90}}
\multiput(0.5,0)(-1.4,0){3}{\color{red}\line(-3,-2){37.5}}
\multiput(-89,-88)(-1.4,0){3}{\color{red}\line(-3,-2){37.5}}
\multiput(-125.5,-25)(-1,0){3}{\color{red}\line(0,-1){90}}
\multiput(0,-89.5)(0,-1){3}{\color{red}\line(-1,0){90}}
\multiput(-37.5,-24.5)(0,1){3}{\color{red}\line(-1,0){90}}
\end{picture}}
\put(255,240){ 
\begin{picture}(270,225)
\multiput(-90,-90)(0,-1.5){60}{\color{green}\circle*{3}}
\multiput(0,-180)(-1.5,0){60}{\color{green}\circle*{3}}
\multiput(0,-180)(-0.625,-0.415){60}{\color{green}\circle*{3}}
\multiput(-37.5,-115)(0,-1.5){60}{\color{green}\circle*{3}}
\multiput(-90,-92)(-0.625,-0.415){60}{\color{green}\circle*{3}}
\multiput(-37.5,-117)(-1.5,0){60}{\color{green}\circle*{3}}
\end{picture}}
\put(255,240){ 
\begin{picture}(270,225)
\multiput(-75,-140)(0,-1.5){60}{\color{blue}{\circle*{3}}} 
\multiput(-75,-230)(-1.5,0){60}{\color{blue}{\circle*{3}}}
\multiput(-37.5,-115)(-0.625,-0.415){60}{\color{blue}{\circle*{3}}}
\multiput(-127.5,-205)(-0.625,-0.415){60}{\color{blue}{\circle*{3}}}
\multiput(-37.5,-113)(-1.5,0){60}{\color{blue}{\circle*{3}}}
\multiput(-125.5,-115)(0,-1.5){60}{\color{blue}{\circle*{3}}}
\end{picture}}
\put(132,125){\circle*{10}}
\end{picture}}
\caption{The chiral polyhedron $P(1,0)$ of type $\{6,6\}$ in the neighborhood of a single vertex.}
\label{poly66}
\end{figure}
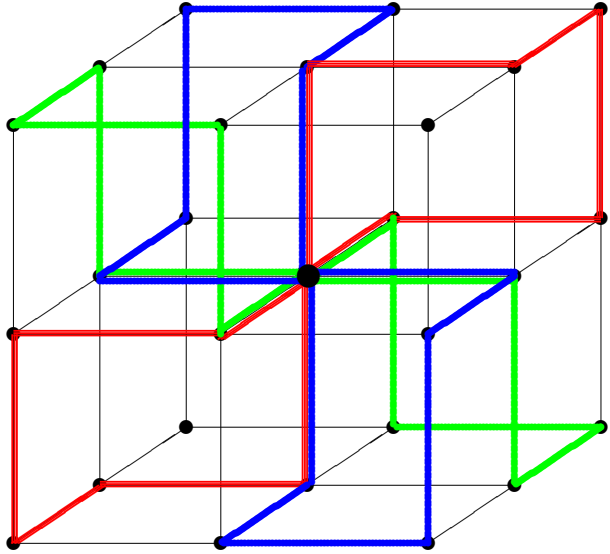

\begin{figure}
\centering
\resizebox{180pt}{255pt}{
\begin{picture}(190,270)
\multiput(10,50)(0,60){3}{
\begin{picture}(230,60)
\thinlines
\multiput(0,0)(60,0){3}{\circle*{1}}
\multiput(25,16.66)(60,0){3}{\circle*{1}}
\multiput(50,33.3)(60,0){3}{\circle*{1}}
\multiput(0,0)(25,16.66){3}{\line(1,0){120}}
\multiput(0,0)(60,0){3}{\line(3,2){50}}
\end{picture}}
\put(10,50){
\begin{picture}(230,100)
\thinlines
\multiput(0,0)(60,0){3}{\line(0,1){120}}
\multiput(25,16.66)(60,0){3}{\line(0,1){120}}
\multiput(50,33.33)(60,0){3}{\line(0,1){120}}
\end{picture}}
\put(10,50){ 
\begin{picture}(260,250)
\thinlines
\put(85,76.66){\vector(1,0){130}}
\put(218,74.66){$\xi_2$}
\put(85,76.66){\line(-1,0){130}}
\put(85,76.66){\vector(0,1){130}}
\put(82,211.5){$\xi_3$}
\put(85,76.66){\line(0,-1){130}}
\put(85,76.66){\line(3,2){55}}
\put(85,76.66){\vector(-3,-2){65}}
\put(10.5,30){$\xi_1$}
\put(85,76.66){\circle*{2}}
\put(205,76.66){\circle*{6}}
\put(-35,76.66){\circle*{6}}
\put(85,196.66){\circle*{6}}
\put(85,-43.33){\circle*{6}}
\put(135,110){\circle*{6}}
\put(35,43.33){\circle*{6}}
\put(0,0){\circle*{6}}
\put(120,0){\circle*{6}}
\put(0,120){\circle*{6}}
\put(50,153.33){\circle*{6}}
\put(170,153.33){\circle*{6}}
\put(170,33.33){\circle*{6}}
\put(85,76.66){\circle*{6}}
\end{picture}}
\put(10,50){ 
\begin{picture}(260,250)
\multiput(0,2)(1.415,1.275){60}{\color{red}\circle*{3}}
\multiput(0,118)(1.7,-0.865){50}{\color{red}\circle*{3}}
\multiput(0,0)(-0.875,1.9){40}{\color{red}\circle*{3}}
\multiput(-35.5,76.66)(1.165,1.44){30}{\color{red}\circle*{3}}
\end{picture}}
\put(10,50){ 
\begin{picture}(260,250)
\multiput(0,121)(1.7,-0.865){50}{\color{green}\circle*{3}}
\multiput(48,153.33)(0.875,-1.94){40}{\color{green}\circle*{3}}
\multiput(0,120)(1.41,1.27){60}{\color{green}\circle*{3}}
\multiput(50,153.33)(1.16,1.44){30}{\color{green}\circle*{3}}
\end{picture}}
\put(10,50){ 
\begin{picture}(260,250)
\multiput(0,-2)(1.415,1.275){60}{\color{blue}\circle*{3}}
\multiput(0,0)(1.165,1.44){30}{\color{blue}\circle*{3}}
\multiput(35,43.33)(1.7,-0.865){50}{\color{blue}\circle*{3}}
\multiput(84,76.66)(0.875,-1.94){40}{\color{blue}\circle*{3}}
\end{picture}}
\put(180,203.33){ 
\begin{picture}(260,250)
\multiput(0,-2)(-1.415,-1.275){60}{\color{red}\circle*{3}}
\multiput(0,-118)(-1.7,+0.865){50}{\color{red}\circle*{3}}
\multiput(0,0)(0.875,-1.9){40}{\color{red}\circle*{3}}
\multiput(35.5,-76.66)(-1.165,-1.44){30}{\color{red}\circle*{3}}
\end{picture}}
\put(180,203.33){ 
\begin{picture}(260,250)
\multiput(0,-121)(-1.7,0.865){50}{\color{green}\circle*{3}}
\multiput(-48,-153.33)(-0.875,1.94){40}{\color{green}\circle*{3}}
\multiput(0,-120)(-1.41,-1.27){60}{\color{green}\circle*{3}}
\multiput(-50,-153.33)(-1.16,-1.44){30}{\color{green}\circle*{3}}
\end{picture}}
\put(180,203.33){ 
\begin{picture}(260,250)
\multiput(0,2)(-1.415,-1.275){60}{\color{blue}\circle*{3}}
\multiput(0,0)(-1.165,-1.44){30}{\color{blue}\circle*{3}}
\multiput(-35,-43.33)(-1.7,0.865){50}{\color{blue}\circle*{3}}
\multiput(-84,-76.66)(-0.875,1.94){40}{\color{blue}\circle*{3}}
\end{picture}}
\put(99,126.66){\circle*{10}}
\end{picture}}
\caption{The chiral polyhedron $Q(1,1)$ of type $\{4,6\}$ in the neighborhood of a single vertex.}
\label{poly46}
\end{figure}
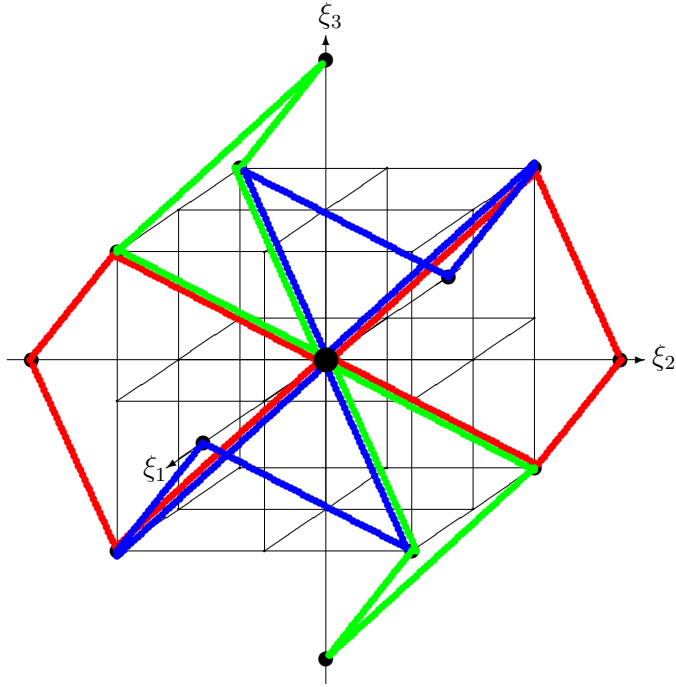

\medskip

\noindent
{\sl Finite-faced Polyhedra\/}
\smallskip

\noindent
In the three families with finite faces, the chiral polyhedra have Schl\"afli symbols $\{6,6\}$, $\{4,6\}$, or $\{6,4\}$, respectively, and both the vertex-figures and the faces are skew regular polygons (hexagons or tetragons). The latter property sets these polyhedra apart from regular polyhedra. In fact, every regular polyhedron with finite faces must necessarily have either planar faces or planar vertex-figures. The two exceptional, regular polyhedra in each of the three extended families result from certain normally excluded parameter combinations and indeed have planar faces or planar vertex-figures.

Figures~\ref{poly66} and~\ref{poly46} show the neighborhood of a single vertex in two chiral polyhedra of types $\{6,6\}$ and $\{4,6\}$, respectively. Both polyhedra have their vertices in the standard cubical lattice of $\E$. In the polyhedron $P(1,0)$ of type $\{6,6\}$ shown in Figure~\ref{poly66}, the faces are Petrie polygons of cubes in the corresponding cubical tessellation; there are six faces meeting at each vertex, with the vertex-figure given by a Petrie polygon of the octahedral vertex-figure of the cubical tessellation at that vertex. The polyhedron $P(1,0)$ is one of only four polyhedra in the extended family of type $\{6,6\}$ whose edges are edges of the cubical tessellation; the others are $P(-1,0)$ and $P(0,\pm 1)$, which are congruent to $P(1,0)$. In the polyhedron $Q(1,1)$ of type $\{4,6\}$ shown in Figure~\ref{poly46}, the edges are main diagonals of cubes in the cubical tessellation, with the four cubes that determine a tetragonal face sharing one edge. 

The finite-faced chiral polyhedra are parametrized up to congruence by a pair of integral parameters $a,b$ (or $c,d$), not both zero, and relatively prime if both non-zero. The distinguished generators $S_1,S_2$ of the symmetry group are rotatory reflections defined in terms of these parameters, and the symmetry $T$ is a half-turn. 

For the polyhedra in the extended family of type $\{6,6\}$, denoted $P(a,b)$, the symmetries $S_1$, $S_2$ and $T$ are given by  
\[ \begin{array}{rccl}
S_{1}\colon & (x_{1},x_{2},x_{3}) & \mapsto & (-x_{2},x_{3},x_{1}) + (0,-b,-a),\\ 
S_{2}\colon & (x_{1},x_{2},x_{3}) & \mapsto & -(x_{3},x_{1},x_{2}),  \\
T\colon & (x_{1},x_{2},x_{3}) & \mapsto &  (-x_{1},x_{2},-x_{3}) + (a,0,b).  
\end{array} \]
The polyhedron $P(a,b)$ itself is obtained from the group $G(a,b)$ generated by $S_1,S_2$ by means of a ``chiral variant" of  Wythoff's construction for regular polyhedra~\cite{crp,arp}. The base vertex is the origin $o$ (fixed by $S_2$); the base edge is given by $\{o,u\}$ with $u:= T(o) = (a,0,b)$; and the base face is determined by the orbit of $o$ under $\langle S_1\rangle$ and has the vertices  
\[ \begin{array}{l}
(0,0,0),(0,-b,-a),(b,-a-b,-a),(a+b,-a-b,-a+b),(a+b,-a,b),(a,0,b), 
\end{array}\]
listed in cyclic order. Thus the integers $a,b$ parametrize the location of the vertex of the base edge distinct from the base vertex. The vertex-figure of $P(a,b)$ at $o$ is determined by the orbit of $u$ under $\langle S_2\rangle$ and has the vertices
\[\begin{array}{l}
(a,0,b),(-b,-a,0),(0,b,a),(-a,0,-b),(b,a,0),(0,-b,-a),
\end{array}\] 
listed in cyclic order. The six faces of $P(a,b)$ containing the base vertex $o$ are given by the orbit of the base face under $\langle S_2\rangle$. Each face is a (generally) skew hexagon with vertices given by one half of the vertices of a hexagonal prism. The vertices, edges, and faces of the entire polyhedron $P(a,b)$ are the elements in the orbits of the base vertex, base edge, or base face under $G(a,b)$. The polyhedron $P(1,0)$ is shown in Figure~\ref{poly66}.

The polyhedra $P(a,b)$ are chiral, unless $b=\pm a$. In the two exceptional cases the polyhedra acquire additional symmetries and become regular; these are the two regular polyhedra in the extended family of type $\{6,6\}$. 
The case $b=a$ gives the Petrie-Coxeter polyhedron $\{6,6\,|\,3\}$, with convex faces and skew vertex-figures, shown in Figure~\ref{petcox}. When $b=-a$ we obtain the pure regular polyhedron $\{6,6\}_4$, with skew faces and convex vertex-figures; its vertices comprise the vertices in one set of alternate vertices of the Petrie-Coxeter polyhedron $\{4,6\,|\,4\}$, while its faces are the vertex-figures at the vertices in the other set of alternate vertices of $\{4,6\,|\,4\}$.

The polyhedra $P(a,b)$ and $P(b,a)$ are congruent. They also are geometric duals of each other, meaning that the face centers of one are the vertices of the other. Thus $P(a,b)$ is geometrically self-dual, in the sense that its dual, $P(a,b)^{*}=P(b,a)$, is congruent to~$P(a,b)$. 

The two extended families of type $\{4,6\}$ or $\{6,4\}$ behave similarly. The corresponding polyhedra are denoted by $Q(c,d)$ and $Q(c,d)^*$, respectively, where the star indicates that the polyhedron $Q(c,d)^*$ of type $\{6,4\}$ is the dual of the polyhedron $Q(c,d)$ of type $\{4,6\}$ (again with the roles of vertices and face centers interchanged). In either case,  $S_1$ and $S_2$ are rotatory reflections accounting for the (generally) skew faces and vertex-figures. The polyhedron $Q(1,1)$ is shown in Figure~\ref{poly46}.

It is quite remarkable that essentially any two finite-faced polyhedra of the same type are non-isomorphic; more precisely, $P(a,b)$ and $P(a',b')$ are combinatorially isomorphic if and only if $(a',b')=\pm (a,b),\pm (b,a)$, and similarly, $Q(c,d)$ and $Q(c',d')$ (and hence their duals, $Q(c,d)^*$ and $Q(c',d')^*$) are combinatorially isomorphic if and only if $(c',d')=\pm (c,d),\pm (-c,d)$. Thus there are very many different isomorphism classes of finite-faced chiral polyhedra. If the enumeration of these polyhedra is carried out by similarity classes, rather than congruence classes, there is just one single rational parameter, namely $a/b$ or $c/d$ (when $b,d\neq 0$). The abundance of isomorphism classes of finite-faced chiral polyhedra then means that the similarity classes exhibit a very strong discontinuity, in that any small change in the rational parameter produces a new similarity class in which the polyhedra are not combinatorially isomorphic to those in the original similarity class. At present no intuitive indicator is known that can detect combinatorial non-isomorphisms among the finite chiral polyhedra.

Moreover, as was shown in \cite{pelwei}, the finite-faced chiral polyhedra are also combinatorially chiral, meaning that the underlying abstract polyhedron is abstractly chiral, in the sense that its combinatorial automorphism group has two flag-orbits such that adjacent flags are in distinct orbits \cite{SW1}.
\medskip

\noindent
{\sl Helix-faced Polyhedra\/}
\smallskip

\noindent
In the three helix-faced families, the polyhedra have Schl\"afli symbols $\{\infty,3\}$, $\{\infty,3\}$, or $\{\infty,4\}$, respectively, and the helical faces are over triangles, squares, and triangles. The vertex-figures are convex. Each family has two real-valued parameters $a,b$ or $b,c$ that cannot both be zero; the polyhedra are denoted by $P_1(a,b)$, $P_2(c,d)$ and $P_3(c,d)$, respectively. Unlike in the case of chiral polyhedra with finite faces (where the parameters were integers), the discreteness does not impose any further restrictions on these parameters. Again, each (extended) family contains two regular polyhedra, but now one polyhedron is finite and is given by a crystallographic Platonic solid, that is, a tetrahedron, cube, or octahedron, respectively,

As before, the symmetry groups are generated by two symmetries $S_1$ and $S_2$ determined by a base flag. Except when a polyhedron is finite, $S_1$ is a screw motion moving the vertices of the helical base face one step along the face; 
this would also extend to the finite polyhedron, if we allowed the translation component to be trivial. Now $S_2$ is an ordinary rotation in an axis through the base vertex, $o$ (say), which accounts for the planarity of the vertex-figure. The symmetry $T$ is again a half-turn about an axis passing through the midpoint of the base edge.

For the polyhedra $P_2(c,d)$ in the second family of type $\{\infty,3\}$, with helical faces over squares, the symmetries $S_1$, $S_2$ and $T$ are given by 
\[ \begin{array}{rccl}
S_{1}\colon & (x_{1},x_{2},x_{3}) & \mapsto & (-x_{3},x_{2},x_{1}) + (d,c,-c),\\ 
S_{2}\colon & (x_{1},x_{2},x_{3}) & \mapsto & (x_{2},x_{3},x_{1}),  \\
T\colon & (x_{1},x_{2},x_{3}) & \mapsto & (x_{2},x_{1},-x_{3}) + (c,-c,d).  
\end{array} \]
The entire polyhedron $P_{2}(c,d)$ is again an orbit structure produced from its symmetry group $\langle S_1,S_2\rangle$ by a variant of Wythoff's construction. Its base vertex is $o$, and the vertex in its base edge distinct from $o$ is given by $T(o) =  (c,-c,d)$. The vertex-set of the base helical face is given by 
\[ \{(c,-c,d),(0,0,0),(d,c,-c),(c+d,2c,-c+d)\} + \mathbb{Z}\!\cdot\! t\]
with $t:=(0,4c,0)$; the notation means that the four points listed on the left side are successive vertices of this face, whose translates by integral multiples of $t$ comprise all the vertices of the face. The screw motion $S_1$ shifts the base face one step along itself, and since $S_{1}^4$ is a translation, the helical faces are over squares (when $c\neq 0$). The vertices, edges, and faces of the entire polyhedron $P_{2}(c,d)$ are the elements in the orbits of the base vertex, base edge, or base face under the symmetry group. The real parameters $c,d$ specify the location of the vertex $T(o)$ of $P_{2}(c,d)$, which lies in the plane $x_{2}=-x_{1}$ and a~priori can be any point in that plane other than the origin.

When $c=0$ we obtain a finite polyhedron $P_{2}(0,d)$, namely, a cube $\{4,3\}$; in fact, in this case $S_{1}^4$ is the identity mapping and the base face is a square. In a sense, a chiral helix-faced polyhedron $P_{2}(c,d)$ can be thought of as unraveling the cube; in fact, there is a covering $P_{2}(c,d)\! \mapsto\! \{4,3\}$ under which each helical face of $P_{2}(c,d)$ is ``compressed" (like a spring) into the square over which it had been rising.

It is very remarkable that, for each helix-faced family, each of its chiral polyhedra is combinatorially isomorphic to the infinite regular polyhedron in this family~\cite{pelwei}.  Thus each helix-faced chiral polyhedron is combinatorially regular (as an abstract polyhedron), but is not geometrically regular. In fact, the chiral polyhedra in each family can be viewed as ``chiral deformations" of the infinite regular polyhedron in this family; at the other extreme they can also be moved into the finite regular polyhedron (Platonic solid) in this family. Figure~\ref{chirdeform} illustrates this behavior for the chiral polyhedra $P_{2}(c,d)$ (with $c\neq 0$), which can be deformed into the helix-faced regular polyhedron $P_{2}(1,0)$ shown in Figure~\ref{manhattan}. The details involve two steps, of which only the second is shown. The first step is simple:\ a polyhedron $P_{2}(c,d)$ is rescaled to obtain a similar polyhedron of the form $P_{1}(1,e)$, with $e=d/c$. In the second step the regular polyhedron $P_{2}(1,0)$ is continuously transformed to the chiral polyhedron $P_{2}(1,e)$ (if $e\neq 0$). This second step, which works for any $e$, changes the vertical (say) stacks of cubes involved in the construction of $P_{2}(1,0)$ in a continuous fashion, so as to provide the transition from $P_{2}(1,0)$ to $P_{2}(1,e)$ (see~\cite{pelwei}). These phenomena are in sharp contrast to those observed for the finite-faced chiral polyhedra, where each chiral polyhedron is also combinatorially chiral. 
\bigskip

\noindent
{\bf Polygonal Complexes\/}
\medskip

\noindent
Polygonal complexes form an even broader class of skeletal structures than polyhedra, and are of natural interest to crystallographers. Like polyhedra they are made up of vertices, edges, and polygon faces; however, unlike in polyhedra, more than two faces can meet at an edge. 

More precisely, following~\cite{pesch1,pesch2}, a polygonal complex $K$ in $\mathbb{E}^3$ consists of a set $V$ of points, called vertices, a set $E$ of line segments, called edges, and a set $F$ of polygons, called faces, with the following properties:\\[.04in]
(a)\ The graph $(V,E)$, called the edge graph of $K$, is connected.\\[.01in]
(b)\ The vertex-figure at each vertex of $K$ is connected. By the vertex-figure of $K$ at a vertex~$v$ we mean the graph, possibly with multiple edges, whose vertices are the neighbors of $v$ in the edge graph of $K$ and whose edges are the line segments $(u,w)$, where $(u,v)$ and $(v,w)$ are edges of a common face of $K$.\\[.01in]
(c)\ Each edge of $K$ is contained in exactly $r$ faces of $K$, for a fixed number $r\geq 2$.\\[.01in]
(d)\ $K$ is discrete, meaning that each compact subset of $\mathbb{E}^3$ meets only finitely many faces of~$K$.

A simple example of a polygonal complex with $r=4$ is illustrated in Figure~\ref{2skeleton}. This is the $2$-skeleton of the standard cubical tessellation, consisting of the vertices, edges, and square faces of the tessellation; each edge lies in four squares. The vertex-figure is the edge graph of the octahedron.

\begin{figure}
\centering
\resizebox{260pt}{225pt}{
\begin{picture}(130,125)
\multiput(7.5,0)(0,30){3}{
\begin{picture}(120,30)
\thicklines
\multiput(0,0)(30,0){3}{\circle*{2}}
\multiput(12.5,8.33)(30,0){3}{\circle*{2}}
\multiput(25,16.66)(30,0){3}{\circle*{2}}
\multiput(0,0)(12.5,8.33){4}{\line(1,0){60}}
\multiput(0,0)(30,0){3}{\line(3,2){37.5}}
\end{picture}}
\put(7.5,0){
\begin{picture}(120,110)
\thicklines
\multiput(0,0)(30,0){3}{\line(0,1){60}}
\multiput(12.5,8.33)(30,0){3}{\line(0,1){60}}
\multiput(25,16.66)(30,0){3}{\line(0,1){60}}
\multiput(37.5,25)(30,0){3}{\line(0,1){60}}
\end{picture}}
\thicklines
\put(53,39.5){\color{red}{\line(3,2){12.5}}}
\put(53,38.75){\color{red}{\line(3,2){12.5}}}
\put(53,38){\color{red}{\line(3,2){12.5}}}
\put(52,39.5){\color{red}{\line(3,2){12.5}}}
\put(52,36.5){\color{red}{\line(3,2){12.5}}}
\put(22,38){\color{red}{\line(1,0){60}}}
\put(22,38){\color{red}{\line(3,2){12.5}}}
\put(82,38){\color{red}{\line(3,2){12.5}}}
\put(34.5,46.33){\color{red}{\line(1,0){60}}}
\put(53,9){\color{red}{\line(0,1){60}}}
\put(22,38){\color{red}{\line(3,2){12.5}}}
\put(52,68){\color{red}{\line(3,2){12.5}}}
\put(52,8){\color{red}{\line(3,2){12.5}}}
\put(65.5,15.33){\color{red}{\line(0,1){60}}}
\end{picture}}
\caption{The $2$-skeleton of the cubical tessellation.}
\label{2skeleton}
\end{figure}
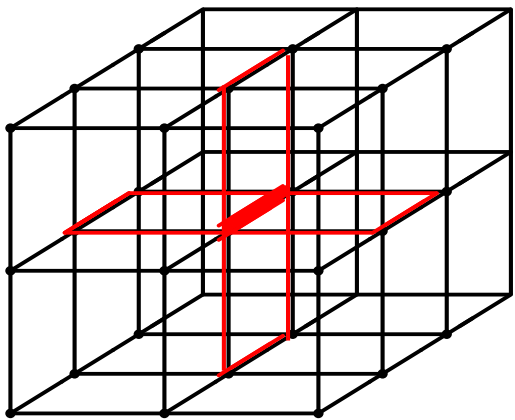

A polyhedron is a polygonal complex with $r=2$. When viewed purely combinatorially, the set of all vertices, edges, and faces of a polygonal complex $K$, ordered by inclusion, is an incidence complex of rank $3$ in the sense of~\cite{ds}; this is an abstract polyhedron if $K$ is a polyhedron. In a polygonal complex $K$, the vertex-figures must necessarily be allowed to have multiple edges, to account for the possibility that adjacent edges of a face of $K$ are adjacent edges of more than one face of $K$. By the discreteness, the vertex-figures are finite graphs. 

A polygonal complex $K$ is regular if its symmetry group $G(K)$ acts transitively on the flags (incident vertex-edge-face triples). As for regular polyhedra, the symmetry group of a regular complex $K$ acts transitively, separately, on the sets of faces, edges, and vertices, so in particular $K$ has mutually congruent faces, edges, and vertex-figures, respectively. The faces of $K$ are finite or infinite regular polygons. 

The symmetry group $G(K)$ of a regular polygonal complex $K$ may or may not be simply flag-transitive (that is, may or may not have trivial flag stabilizers). In either case $G(K)$ is generated by three subgroups $G_0$, $G_1$ and $G_2$ which are defined by their action on a base flag of $K$. When $K$ is a polyhedron these subgroups have order $2$ and are generated by the reflections $R_0$, $R_1$ and $R_2$ described earlier. In the general case, for $i=0,1,2$, the subgroup $G_i$ is the stabilizer of all elements in the base flag, except for the element of rank $i$ (the base vertex, base edge, or base face, respectively). For example, $G_2$ is the stabilizer of the base vertex and the base edge, and thus keeps the line through the base edge pointwise fixed and acts on a plane perpendicular to this line as a cyclic or dihedral group. 

If the symmetry group $G(K)$ is not simply flag-transitive, then $K$ has planar faces and the base flag stabilizer in $G(K)$ has order $2$ and is generated by the reflection in the plane containing the base face of $K$.  In this case we say that $K$ is non-simply flag-transitive, and that $K$ has face mirrors, meaning that each plane through a face of $K$ is the mirror (fixed point set) of a reflection symmetry of $K$. We call $K$ simply flag-transitive if $G(K)$ is simply flag-transitive.
\medskip

\noindent
{\sl Non-Simply Flag-Transitive Complexes\/}
\smallskip

\noindent
There are just four regular polygonal complex $K$ in $\E$ which are non-simply flag transitive. The existence of face mirrors in this case allows us to recognize $K$ as the $2$-skeleton of a certain type of incidence structure of rank $4$ in $\E$, called a regular $4$-apeirotope~\cite[Ch. 7F]{arp}. The complex $K$ then consists of the vertices, edges, and faces of rank 2 of this $4$-apeirotope. There are eight regular $4$-apeirotopes in $\E$, among them the standard cubical tessellation. The $2$-skeleton of the cubical tessellation shown in Figure~\ref{2skeleton} is the simplest example of a non-simply flag transitive polygonal complex. 

The eight regular 4-apeirotopes fall into  four pairs of ``Petrie-duals", with each pair sharing the same 2-skeleton and producing the same polygonal complex. Thus up to similarity there are four non-simply flag-transitive complexes in $\E$, each with the same symmetry group as its two respective regular $4$-apeirotopes. The $2$-skeleton of the cubical tessellation has square faces and is the only example with finite faces. The three other complexes all have (planar) zigzag faces, with either three or four faces meeting at each edge. Thus non-simply flag-transitive regular polygonal complexes in $\E$ are $2$-skeletons of regular structures of rank $4$ in $\E$.
\medskip

\noindent
{\sl Simply Flag-transitive Complexes\/}
\smallskip

\noindent
The classification of the simply flag-transitive regular polygonal complexes is more involved. In addition to the forty-eight  regular polyhedra there are twenty-one simply flag-transitive complexes in $\mathbb{E}^3$ which are not polyhedra, up to similarity. 

In a simply flag-transitive complex $K$ which is not a polyhedron, there are at least three faces meeting at each edge. The generating subgroups $G_0,G_1,G_2$ of $G(K)$ then take a very specific form:\ $G_0$ is generated by a point, line, or plane reflection $R_0$; the subgroup~$G_1$ is generated by a line or plane reflection $R_1$; and $G_2$ is a cyclic group $C_r$ or dihedral group $D_{r/2}$ (of order $r$), where $r$ denotes the number of faces meeting at an edge of $K$. Again it is convenient to introduce a mirror vector, now of the form $(d_0,d_1)$, recording the dimensions $d_0$ and $d_1$ of the mirrors of $R_0$ and $R_1$, respectively. 

The twenty-one simply flag-transitive complexes are summarized in Table~\ref{tabsimply}. They are denoted $\K_i(j,k)$, where $(j,k)$ indicates the mirror vector and $i$ is the serial number in the list of regular complexes with mirror vector $(j,k)$. There are columns for the pointwise edge stabilizer $G_2$, the number $r$ of faces at each edge, the types of faces and vertex-figures, the vertex-set, and the corresponding net (which is described below). 

\begin{table}
\centering
\begin{tabular}{|l|l|l|l|l|l|l|}
\hline
Complex&$G_2$ & $r$ &Face&Vertex-Figure&Vertex-Set& Net\\
\hline\hline
$K_1(1,2)$& $D_2$  & $4$ &$4_s$ & cuboctahedron&$\Lambda_{2}$&$\!{\rm fcu}$\\
$K_2(1,2)$& $C_3$& $3$ &$4_s$ & cube&$\Lambda_{3}$&$\!{\rm bcu}$\\
$K_3(1,2)$& $D_3$& $6$   &$4_s$ &double cube&$\Lambda_{3}$&$\!{\rm bcu}$\\
$K_4(1,2)$& $D_2$& $4$ & $6_s$&octahedron&$\Lambda_{1}$&$\!{\rm pcu}$\\
$K_5(1,2)$& $D_2$& $4$ &$6_s$&double square&$V$&$\!{\rm nbo}$\\
$K_6(1,2)$& $D_4$& $8$ &$6_s$&double octahedron&$\Lambda_{1}$&$\!{\rm pcu}$\\
$K_7(1,2)$& $D_3$& $6$ &$6_s$&double tetrahedron&$W$&$\!{\rm dia}$\\
$K_8(1,2)$& $D_2$& $4$ &$6_s$&cuboctahedron&$\Lambda_{2}$&$\!{\rm fcu}$\\
\hline
$K_1(1,1)$& $D_3$& $6$ &$\infty_3$&double cube&$\Lambda_{3}$&$\!{\rm bcu}$\\
$K_2(1,1)$& $D_2$& $4$ &$\infty_3$&double square&$V$&$\!{\rm nbo}$\\
$K_3(1,1)$& $D_4$& $8$ &$\infty_3$&double octahedron&$\Lambda_{1}$ & $\!{\rm pcu}$ \\
$K_4(1,1)$& $D_3$ & $6$ &$\infty_4$& double tetrahedron&$W$ &$\!{\rm dia}$\\
$K_5(1,1)$& $D_2$& $4$ &$\infty_4$&ns-cuboctahedron&$\Lambda_{2}$&$\!{\rm fcu}$ \\
$K_6(1,1)$& $C_3$& $3$ &$\infty_4$&tetrahedron& $W$&$\!{\rm dia}$\\
$K_7(1,1)$& $C_4$& $4$ &$\infty_3$&octahedron&$\Lambda_{1}$ &$\!{\rm pcu}$\\
$K_8(1,1)$& $D_2$& $4$ &$\infty_3$&ns-cuboctahedron&$\Lambda_{2}$&$\!{\rm fcu}$ \\
$K_9(1,1)$& $C_3$& $3$ &$\infty_3$&cube&$\Lambda_{3}$&$\!{\rm bcu}$\\
\hline
$K(0,1)$& $D_2$& $4$ &$\infty_2$&ns-cuboctahedron&$\Lambda_{2}$&$\!{\rm fcu}$ \\
\hline
$K(0,2)$& $D_2$& $4$ &$\infty_2$&cuboctahedron& $\Lambda_{2}$&$\!{\rm fcu}$ \\
\hline
$K(2,1)$& $D_2$& $4$ & $6_c$ &ns-cuboctahedron&$\Lambda_{2}$&$\!{\rm fcu}$\\
\hline
$K(2,2)$& $D_2$& $4$ &$3_c$&cuboctahedron&$\Lambda_{2}$&$\!{\rm fcu}$\\
\hline
\end{tabular}
\caption{The twenty-one simply flag-transitive regular complexes in $\mathbb{E}^3$ which are not polyhedra.}
\label{tabsimply}
\end{table}

In the face column we use symbols like $p_c$, $p_s$, $\infty_2$, or $\infty_k$ with $k=3$ or $4$, respectively, to indicate that the faces are \underbar{c}onvex $p$-gons, \underbar{s}kew $p$-gons, planar zigzags, or helical polygons over $k$-gons. (A planar zigzag is viewed as a helical polygon over a $2$-gon.)  An entry in the vertex-figure column listing a solid figure in $\E$ is meant to represent the geometric edge-graph of this figure, with ``double" indicating the double edge-graph (the edges have multiplicity 2). An entry ``ns-cuboctahedron" stands for the edge graph of a ``\underbar{n}on-\underbar{s}tandard cuboctahedron", a certain realization of the ordinary cuboctahedron with equilateral triangles and skew square faces (see \cite[p.~2041]{pesch2}).  

For all but five complexes the vertex-set is a lattice, namely one of the following:\ the standard cubic lattice $\Lambda_{1}:=\mathbb{Z}^{3}$; the face-centered cubic lattice $\Lambda_{2}$, with basis $(1,1,0)$, $(-1,1,0)$, $(0,-1,1)$, consisting of all integral vectors with even coordinate sum; or the body-centered cubic lattice $\Lambda_3$, with basis $(2,0,0)$, $(0,2,0)$, $(1,1,1)$. The vertex-sets of the five exceptional complexes are either  
\[ V:=\Lambda_{1}\!\setminus\! ((0,0,1)\!+\!\Lambda_{3})\]
or
\[ W:= 2\Lambda_{2} \cup ((1,-1,1)\!+\!2\Lambda_{2}). \]

\begin{figure}
\centering
\resizebox{250pt}{206pt}{
\begin{picture}(130,110)
\multiput(10,0)(0,40){3}{
\begin{picture}(100,25)
\thinlines
\multiput(0,0)(40,0){3}{\circle*{2}}
\multiput(18,12)(40,0){3}{\circle*{2}}
\multiput(36,24)(40,0){3}{\circle*{2}}
\multiput(0,0)(18,12){3}{\line(1,0){80}}
\multiput(0,0)(40,0){3}{\line(3,2){36}}
\end{picture}}
\put(10,0){
\begin{picture}(180,140)
\thinlines
\multiput(0,0)(40,0){3}{\line(0,1){80}}
\multiput(18,12)(40,0){3}{\line(0,1){80}}
\multiput(36,24)(40,0){3}{\line(0,1){80}}
\end{picture}}
\put(10,0){
\begin{picture}(100,80)
\thicklines
\put(0,40){\color{green}\line(1,1){40}}
\multiput(40,81)(.45,-.7){40}{\color{green}\circle*{1}}          
\put(58,52){\color{green}\line(-1,1){40}}
\multiput(18,93.5)(-.3,-.866){60}{\color{green}\circle*{1}}
\end{picture}}
\put(10,0){
\begin{picture}(100,80)
\thicklines
\multiput(0,40)(0.98,.2){60}{\color{yellow}\circle*{1}} 
\multiput(40,79)(.45,-.7){40}{\color{yellow}\circle*{1}} 
\multiput(40,80)(-.55,.3){40}{\color{yellow}\circle*{1}}         
\multiput(18,90)(-.3,-.866){60}{\color{yellow}\circle*{1}}
\end{picture}}
\put(10,0){
\begin{picture}(100,80)
\thicklines
\put(40,80){\color{red}\line(1,-1){40}}
\multiput(80,42)(.3,.866){60}{\color{red}\circle*{1}}
\put(58,52){\color{red}\line(1,1){40}}
\multiput(40,82.5)(.45,-.68){42}{\color{red}\circle*{1}}  
\end{picture}}
\put(10,0){
\begin{picture}(100,80)
\thicklines
\multiput(80,40)(-.55,.3){40}{\color{blue}\circle*{1}}      
\multiput(80,38)(.3,.866){60}{\color{blue}\circle*{1}}
\multiput(40,80)(0.98,.2){60}{\color{blue}\circle*{1}} 
\multiput(39,77.5)(.45,-.68){40}{\color{blue}\circle*{1}}  
\end{picture}}
\put(10,0){
\begin{picture}(80,80)
\put(0,40){\circle*{5}}
\put(18,92){\circle*{5}}
\put(80,40){\circle*{5}}
\put(98,92){\circle*{5}}
\put(40,80){\circle*{5}}
\put(58,52){\circle*{5}}
\end{picture}}
\end{picture}}
\caption{The four tetragon faces of $K_{1}(1,2)$ sharing an edge.}
\label{figk112}
\end{figure}
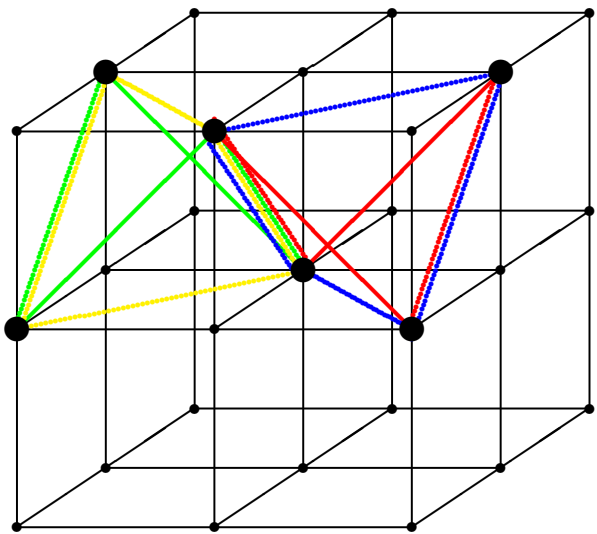

Figure~\ref{figk112} illustrates the complex $K_{1}(1,2)$ with tetragons as faces, four meeting at an edge. The faces are the Petrie polygons of regular tetrahedra inscribed in the cubes of the cubical tessellation of $\E$, where tetrahedra in adjacent cubes are images of each other under the reflection in the common square face of the cubes. The vertex-set of $K_{1}(1,2)$ is the face-centered cubic lattice and the vertex-figure is the edge graph of the ordinary cuboctahedron. 

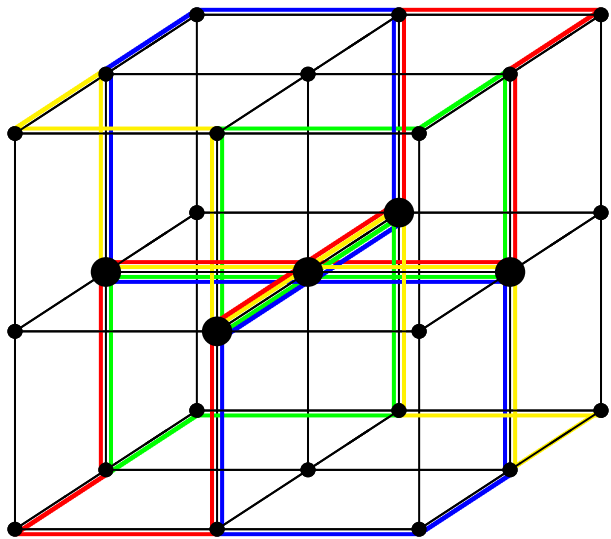
\begin{figure}
\centering
\resizebox{275pt}{206pt}{
\begin{picture}(140,110)
\multiput(10,0)(0,40){3}{
\begin{picture}(100,25)
\thinlines
\multiput(0,0)(40,0){3}{\circle*{3}}
\multiput(18,12)(40,0){3}{\circle*{3}}
\multiput(36,24)(40,0){3}{\circle*{3}}
\multiput(0,0)(18,12){3}{\line(1,0){80}}
\multiput(0,0)(40,0){3}{\line(3,2){36}}
\end{picture}}
\put(10,0){
\begin{picture}(180,140)
\thinlines
\multiput(0,0)(40,0){3}{\line(0,1){80}}
\multiput(18,12)(40,0){3}{\line(0,1){80}}
\multiput(36,24)(40,0){3}{\line(0,1){80}}
\end{picture}}
\put(10,0){
\begin{picture}(100,80)
\thicklines
\put(41,40){\color{green}\line(0,1){40}}
\put(40,81){\color{green}\line(1,0){40}}
\multiput(80,81)(.45,.3){40}{\color{green}\circle*{1}}          
\put(97,52){\color{green}\line(0,1){40}}
\put(98,51){\color{green}\line(-1,0){40}}
\multiput(58,51)(-.45,-.3){40}{\color{green}\circle*{1}}
\end{picture}}
\put(10,0){
\begin{picture}(100,80)
\thicklines
\put(58,54){\color{red}\line(1,0){40}}
\put(99,52){\color{red}\line(0,1){40}}
\multiput(98,93)(.45,.3){40}{\color{red}\circle*{1}}    
\put(116,105){\color{red}\line(-1,0){40}}
\put(77,104){\color{red}\line(0,-1){40}}
\multiput(58,54.2)(.45,.3){40}{\color{red}\circle*{1}}
\put(58,50){\color{blue}\line(-1,0){40}} 
\put(19,52){\color{blue}\line(0,1){40}}
\multiput(18,93)(.45,.3){40}{\color{blue}\circle*{1}}
\put(36,105){\color{blue}\line(1,0){40}}
\put(75,104){\color{blue}\line(0,-1){40}}
\multiput(76,61.6)(-.45,-.3){40}{\color{blue}\circle*{1}}
\multiput(58,51)(.45,.3){40}{\color{green}\circle*{1}}
\put(75,64){\color{green}\line(0,-1){40}}
\put(76,23){\color{green}\line(-1,0){40}}
\multiput(36,23)(-.45,-.3){40}{\color{green}\circle*{1}}
\put(17,51){\color{green}\line(1,0){40}}
\put(19,52){\color{green}\line(0,-1){40}}
\multiput(58,53)(-.45,-.3){40}{\color{yellow}\circle*{1}}
\put(58,53){\color{yellow}\line(-1,0){40}}
\put(39,40){\color{yellow}\line(0,1){40}}
\put(40,81){\color{yellow}\line(-1,0){40}}
\multiput(0,81)(.45,.3){40}{\color{yellow}\circle*{1}}
\put(17,52){\color{yellow}\line(0,1){40}}
\put(58,54){\color{red}\line(-1,0){40}}
\put(17,52){\color{red}\line(0,-1){40}}
\multiput(18,11)(-.45,-.3){40}{\color{red}\circle*{1}}    
\put(0,-1){\color{red}\line(1,0){40}}
\put(39,0){\color{red}\line(0,1){40}}
\multiput(40,42)(.45,.3){40}{\color{red}\circle*{1}}
\put(58,50){\color{blue}\line(1,0){40}} 
\put(97,52){\color{blue}\line(0,-1){40}}
\multiput(98,11)(-.45,-.3){40}{\color{blue}\circle*{1}}
\put(80,-1){\color{blue}\line(-1,0){40}}
\put(41,0){\color{blue}\line(0,1){40}}
\multiput(40,38)(.45,.3){40}{\color{blue}\circle*{1}}
\multiput(58,53)(.45,.3){40}{\color{yellow}\circle*{1}}
\put(77,64){\color{yellow}\line(0,-1){40}}
\put(76,23){\color{yellow}\line(1,0){40}}
\put(99,12){\color{yellow}\line(0,1){40}}
\multiput(116,24)(-.45,-.3){40}{\color{yellow}\circle*{1}}
\put(58,53){\color{yellow}\line(1,0){40}}
\put(58,52){\circle*{6}}
\put(18,52){\circle*{6}}
\put(98,52){\circle*{6}}
\put(40,40){\circle*{6}}
\put(76,64){\circle*{6}}
\end{picture}}
\multiput(10,0)(0,40){3}{
\begin{picture}(100,25)
\thinlines
\multiput(0,0)(40,0){3}{\circle*{3}}
\multiput(18,12)(40,0){3}{\circle*{3}}
\multiput(36,24)(40,0){3}{\circle*{3}}
\multiput(0,0)(18,12){3}{\line(1,0){80}}
\multiput(0,0)(40,0){3}{\line(3,2){36}}
\end{picture}}
\end{picture}}
\caption{The eight faces of ${K}_{5}(1,2)$ meeting at a vertex. }
\label{figk512}
\end{figure}

Figure~\ref{figk512} shows the local structure of the polygonal complex $K_{5}(1,2)$. Its faces are Petrie polygons of cubes in the cubical tessellation, four meeting at each edge. Each cube of the tessellation contributes only one of its Petrie polygons as a face. Shown are the eight faces of ${K}_{5}(1,2)$ that have the central vertex in common. The vertex-figure at this vertex is the double square spanned by the four outer black nodes in the horizontal plane through the center. 

There other ways in which Petrie polygons of cubes in the cubical tessellation can assemble to give the faces of a regular polygonal complex. An example, namely the complex ${K}_{4}(1,2)$, was shown in Figure~\ref{figk412}. Its faces are the Petrie polygons of alternate cubes in the cubical tessellation of $\E$; for every cube occupied, all its Petrie polygons occur as faces. The figure only depicts the twelve faces that have the central vertex in common. The vertex-figure of $K_{4}(1,2)$ at this vertex is the edge graph of the octahedron spanned by the six outer black nodes. 
\medskip

Thus, in summary, there are a total of seventy-three regular polygonal complexes in $\E$, namely forty-eight regular polyhedra (the Gr\"unbaum-Dress polyhedra), four non-simply flag-transitive polygonal complexes (arising as $2$-skeletons of regular apeirotopes of rank 4 in $\E$), and twenty-one simply flag-transitive polygonal complexes.
\bigskip

\noindent
{\bf Nets\/}
\medskip

\noindent
In crystal chemistry there is considerable interest in the relationship between crystal nets and regular figures in space~\cite{okeef,okehyd,wells}. Nets are 3-periodic graphs in $\E$ that represent crystal structures, in the simplest form with vertices corresponding to atoms and edges to bonds. There are extensive databases for crystal nets (RCSR, TOPOS) which in particular include the most symmetric examples~\cite{bla,oketal}. 

Formally, a net $N$ in $\E$ is a 3-periodic connected (simple) graph in $\E$, with straight edges. Here $3$-periodicity means that the translation subgroup of the symmetry group $G(N)$ of N is generated by three translations in independent directions. In the present context all nets are uninodal, meaning that $G(N)$ is transitive on the vertices (nodes) of $N$. For a vertex $v$ of $N$, the convex hull of the neighbors of $v$ in $N$ is called the coordination figure of $N$ at $v$. In a highly symmetric net, the coordination figures are highly symmetric convex polyhedra or polygons.

The edge graphs (1-skeletons) of almost all regular polygonal complexes in $\E$ are highly-symmetric nets, with the only exceptions arising from the polyhedra which are not 3-periodic (that is, from the finite polyhedra and planar polyhedra). 
The last column of Table~\ref{tabsimply} lists the nets of the twenty-one simply-flag transitive regular polygonal complexes which are not polyhedra. Only five nets occur, and each of these five occurs more than once. The names {\sl pcu\/}, {\sl fcu\/}, and {\sl bcu\/} stand for the ``primitive cubic lattice" (the standard cubic lattice), the ``face-centered cubic lattice", and the ``body-centered" cubic lattice in $\E$, respectively. The net {\sl dia\/} is the famous diamond net, which is the net of the diamond form of carbon. The notation {\sl nbo\/} signifies the net of Niobium Monoxide, NbO. When any of these five nets is realized as the edge graph of a regular polygonal complex, its coordination figure occurs as the convex hull of the vertex-figure of the complex (ignoring the multiplicity of edges in the vertex-figure). 

\begin{figure}
\centering
\includegraphics[width=3.2in,height=2.88in]{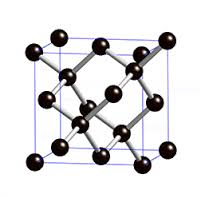}
\caption{The diamond net, occurring as the edge graph of $K_{7}(1,2)$.}
\label{figk712}
\end{figure}

The edge-graph of $K_{7}(1,2)$ is the diamond net, {\sl dia\/}, modeling the diamond crystal shown in Figure~\ref{figk712}~\cite{wiki,wolfram}. The carbon atoms sit at the vertices, and the bonds between adjacent atoms are represented by edges. The ``hexagonal rings" (fundamental circuits) of the net are the faces of $K_{7}(1,2)$.

\end{document}